\documentclass[12pt]{article}
\usepackage{latexsym}                 
\usepackage[T1]{fontenc}
\usepackage{lmodern}
\usepackage[latin1]{inputenc}
\usepackage{latexsym}                 
\usepackage{color}                 
\usepackage{epsfig}
\usepackage{amssymb}
\usepackage{amsthm}
\usepackage{graphicx}
\usepackage{hyperref}
\usepackage{verbatim}
\usepackage{mathptmx}      
\newcommand{\eq}{\end{equation}}
\def\fa{\hbox{ for all }}

\def\b1{\mathbf 1}

\newcommand{\R}{\ensuremath{\mathbb{R}}}

\def\bx{\mathbf{x}}	
\def\boy{\mathbf{y}}

\def\biglf{\par\bigskip\noindent}
\def\RSlabel#1{\label{#1}%
}
\def\RScite#1{\cite{#1}%
}

\newcommand{\bql}[1]{%
 \hfill {\tt ({#1})}
\begin{equation}\RSlabel{#1}%
}
\def\filename#1{}
\begin{document}
\begin{center}
{\Large\bf Adaptive Approximation of Functions with Discontinuities}\\
~\\
Licia Lenarduzzi and Robert Schaback\\
~\\
Version of Nov. 09, 2015
\end{center} 
{\bf Abstract}: 
One of the basic principles of Approximation Theory is that
the quality of approximations 
increase with the smoothness of the function to be approximated.
Functions that are smooth in certain subdomains will have 
good approximations in those subdomains, and these {\em sub-approximations} 
can possibly be calculated efficiently 
in parallel, as long as the subdomains do not
overlap. This paper proposes a class of algorithms that first
calculate sub-approximations on non-overlapping subdomains, 
then extend the subdomains as much as possible 
and finally produce a global solution on the given domain 
by letting the subdomains fill the whole domain. 
Consequently, there will be no Gibbs phenomenon along the 
boundaries of the subdomains. Throughout, the algorithm works for 
fixed scattered input data of the function itself,
not on spectral data,  and it does not resample. 




\section{Introduction}\RSlabel{SecIntro}
Assume that a large set $\{(\bx_i,f_i),i=1,\ldots,N\}$ 
of data is given, where the 
points $\bx_i$ are scattered in $\R^d$ 
and form a set $X$.
We want to find a function 
$u$ that recovers the data on a domain $\Omega$ containing the points, i.e.
$$
\begin{array}{rcl}
u&:& \Omega\to \R,\\
u(\bx_i)&\approx& f_i,\;i=1,\ldots,N.
\end{array} 
$$
We are particularly interested in situations where the data have smooth
interpolants in certain non-overlapping
subdomains $\Omega_j$, 
but not globally. The reason may be 
that there are discontinuities in the  function itself or its
derivatives. Thus a major goal is to identify subdomains
$\Omega_j\subseteq\Omega,\;1\leq j\leq J$ and smooth functions $u_j,\;1\leq
j\leq J$ such that  
$$
\begin{array}{rcl}
u_j&:& \Omega_j\to \R,\\
u_j(\bx_i)&\approx& f_i \fa \bx_i\in X\cap \Omega_j.
\end{array} 
$$
The solution to the  problem  is piecewise defined as
$$
u(\bx):=u_j(\bx)\fa \bx\in \Omega_j,\;1\leq j\leq J.
$$
Our motivation is the well-known fact that errors and convergence rates in
Approximation Theory always improve with increasing smoothness. Thus on each
subdomain we expect to get rather small errors, 
much smaller than if the problem
would have been treated globally, 
where the non-smoothness is a serious limiting effect. 

From the viewpoint of Machine Learning 
\RScite{cristianini-shawe-taylor:2000-1,schoelkopf-smola:2002-1,%
shawetaylor-cristianini:2004-1}
this is a mixture of classification and regression. The domain points have to
be classified in such a way that on each class there is a good regression 
model. The given training data are used for both classification and regression,
but in this case the classification is dependent on the regression,
and the regression is dependent on the classification.

Furthermore, there is a serious amount of geometry hidden behind the problem.
The subdomains should be connected, their interiors should be disjoint,
and the union of their closures should fill the domain completely.
This is why a black-box machine learning approach is not pursued here.
Instead, Geometry and Approximation Theory play a dominant part.
For the same reason, we avoid to calculate edges or fault lines
first, followed by local approximations later. The approximation properties
should determine  the domains and their boundaries, not the other way round.

In particular, {\em localized approximation} will combine 
Geometry and Approximation Theory 
and provide a central tool, together with {\em adaptivity}.
The basic idea is that in the interior of each subdomain,
far away from its boundary, there should be a good 
and simple approximation to the data at each data point 
from the data of its neighbors. 
\section{An Adaptive Algorithm}\RSlabel{SecA3}
Localized approximation will be used 
as the first  
phase of an {\em adaptive algorithm}, constructing
disjoint localized subsets of the data that allow 
good and simple {\em local approximations}. Thus this
``localization'' phase
produces a subset $X^g \subseteq X$ of ``good'' data points
that is the union of disjoint sets $X_1^g,\ldots,X_J^g$  
consisting of data points that allow good approximations 
$u_j^g\in U,\;1\leq j\leq J$
using only the data points in $X_j^g$. In some sense, this is a rough
classification already, but only of data points. 

The goal of the second 
phase is to reduce the number of
unclassified points by enlarging the sets of classified points. 
It is tacitly assumed that the final number of subdomains is already 
obtained by the number $J$ of classes of ``good'' points 
after the first phase. The ``blow--up'' of the sets $X_j^g$ 
should maintain {\em locality} by adding neighboring data points first,
and adding them only if the local approximation $u_j^g$ does not lose 
too much quality after adding that point and changing the approximation.

The second phase usually leaves a small number of ``unsure'' points 
that could not be clearly classified by blowing up the classified sets.
While the blow-up phase focuses on  each single 
set $X_j^g$ in turn and tries to extend it by
looking at all ``unsure'' points for good extension candidates,
the third phase works the other way round. It focuses on each single 
``unsure'' point $\bx_i$ in turn and looks at all sets $X_j^g$ and the local
approximations $u_j$ on these, and assigns the point $\bx_j$ 
to one of the sets  $X_j^g$ so that $u_j(\bx_i)$ is closest to $f(\bx_i)$.  
It is a ``final assignment'' phase that should classify all data
points and it should produce the final sets $X^f_j\supseteq
X^g_j$ of data points. 
The sets $X^f_j$ should be 
disjoint and their union should be $X$.

After phase 3, each 
local approximation $u_j^f\in U$ 
is based on the points in $X_j^f$ only, but there still are no 
well-defined subdomains 
$\Omega_j\supseteq X_j^f$ as domains of $u_j^f$. 
Thus the determination of subdomain boundaries 
from a classification of data points could be 
the task of a fourth phase.
It could, for instance, be handled by any machine 
learning program that uses the classification
as training data and classifies each given point $\bx$ accordingly.
But this paper does not implement a fourth phase,
being satisfied if each approximation $u_j^f$ is good on each set $X_j^f$,
and much better than any global approximation $u^*\in U$ to all data. 

\section{Implementation}\RSlabel{SecIm}
The above description of a three-phase algorithm allows 
a large variation of different implementations that 
compete for efficiency and accuracy. We shall describe 
a basic implementation together with certain minor variants,
and provide numerical examples demonstrating that 
the overall strategy works fine. 

We work on the unit square of $\R^2$ for simplicity
and take a trial space $U$ 
spanned by translates of a fixed positive definite radial kernel $K$.
In our examples, $K$ may be a Gaussian or an inverse multiquadric.
For details on kernels, readers are referred to standard texts
\RScite{buhmann:2003-1,wendland:2005-1,%
schaback-wendland:2006-1,fasshauer-mccourt:2015-1}, for example.
When working on finite subsets of data points, we shall only use
the translates with respect to this subset. Since the kernel $K$
is fixed, also the Hilbert space $H$ is fixed in which
the kernel is reproducing, and we can evaluate the norm 
$\|.\|_K$ of trial functions
cheaply and exactly. 

To implement locality, we 
assume that we have a computationally cheap method that
allows to calculate for each $\bx\in\R^2$ its $n$ nearest neighbors
from $X$. This can, for instance, be done via a range query
after an initialization of a kd-tree data structure \RScite{bentley:1075-1}.  
\subsection{Phase 1: Localization}\RSlabel{SecP1Loc}
This is carried out by a first step picking all data points
with good localized approximation properties, followed by a second step
splitting the set $X^{g}$ of good points into $J$ disjoint sets $X_j^{g}$. 
\subsubsection{Good Data Points}\RSlabel{SecP1GDP}
We assume  that the global fill distance
$$
h(X, \Omega):=\sup_{\boy\in\Omega}
\min_{\bx_k\in X}\|\boy-\bx_k\|_2
$$
of the full set of data points with respect to the full domain $\Omega$
is roughly the same as the local fill distances
$h(X_j^{f},\Omega_j)$ of the final splitting.

The basic idea is to loop over all $N$ data points  of $X$ and 
to calculate for each data point $\bx_i,\;1\leq i\leq N$ a number $\sigma_i$
that is a reliable indicator 
for the quality of {\em localized approximation}. Using a threshold 
$\sigma$, this allows to determine the set $X^{g}\subseteq X$ of ``good''
data points, without splitting it into subsets.  

There are many ways to do
this. 
The implementation of this paper fixes 
a number $n$ of neighbors and loops over all $N$ data points
to calculate for each data point $\bx_i,\;1\leq i\leq N$
\begin{enumerate}
\item the set $N_i$ of their $n$ nearest neighbors from $X$,
\item the kernel-based interpolant $s_i$ of 
the data $(\bx_k, f(x_k))$ for all $n$ neighboring 
data points $\bx_k\in N_i$,
\item the norm $\sigma_i:=\|s_i\|_K$.
\end{enumerate}  
This loop can be executed with roughly ${\cal O}(Nn^3)$ complexity
and ${\cal O}(N+n^3)$ storage, and with easy parallelization, if necessary
at all. A similar indicator would be the error obtained when
predicting $f(\bx_i)$ from the values at the $n$ neighboring points. 

Practical experience shows that the numbers $\sigma_i$ are good indicators of
locality, because adding outliers to a good interpolant usually
increases the error norm dramatically. Many of the $\sigma_i$ 
can be expected to be small, and thus the threshold
$$
\sigma_i< 2 M_\sigma  
$$ 
will be used to determine ``good'' points within the next
splitting step, see Section \ref{SecP1S}, 
where $M_\sigma$ is the median of all $\sigma_i$.
This is illustrated  for a data set 
by  Figure \ref{loglogsigma}:  
it represents, in base  loglog scale, 
the sorted $\{{\bf \sigma}_i\}$ and the 
constant line relevant to the value of the threshold.
\begin{figure}[hbt]
\psfig{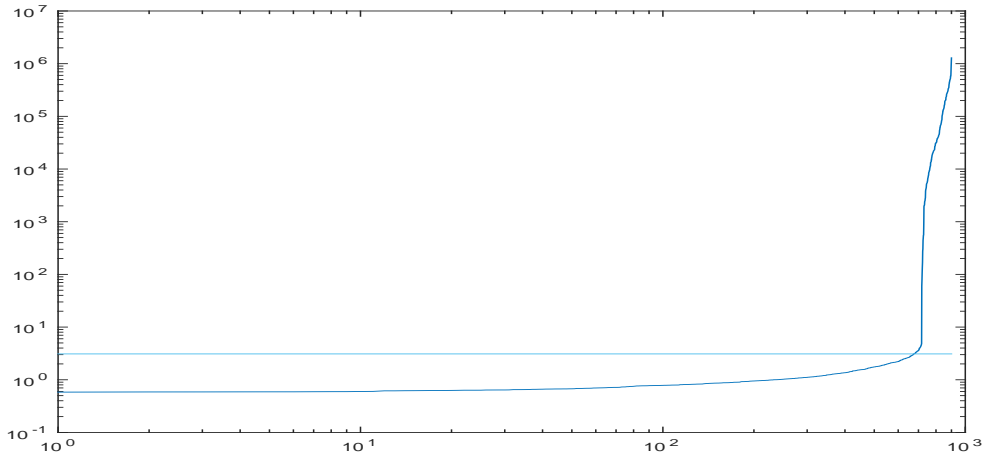}\\
\caption{Loglog: sorted  $\{{\bf \sigma}_i\}$ and  threshold}
\RSlabel{loglogsigma}
\end{figure}

\subsubsection{Splitting}\RSlabel{SecP1S}
The set $X^{g}$ of points with good localization must now be split 
into $J$ disjoint subsets of points that are close to each other. 

We assume that the inner boundaries of the subdomains are 
everywhere clearly determined by 
large values of ${\bf \sigma}_i$.

The implementation of this paper accomplishes the splitting by a variation of
Kruskal's algorithm \RScite{kruskal-1956-1} for calculating minimal spanning trees in
graphs.

The Kruskal algorithm sorts the
edges by increasing weight and starts with an output
 graph that has no edges and no vertices.
When running, it keeps a number of disconnected graphs as the output graph.
It gradually 
adds edges with increasing weight that either connect two previously
disconnected graphs or add an edge to an existing component 
or define a new connected component by that single edge.

In the current implementation the 
edges that connect each $x_i$ with its $n-1$ nearest 
neighbors are collected in an edge list. The edge list
is sorted by
increasing length of the edges 
and then, by $n\mid X \mid$ comparisons, many 
repetitions of edges are removed, and these are all 
repetitions if any two different edges have
different length.

Then the thresholding of the 
$\{{\bf \sigma}_i\}$ by  \[{\bf \sigma}_i< 2 M_{\bf \sigma}\]
is executed, and it is known which points are good and which are bad.

All edges with one or two bad end points are removed from the edge list
with  cost  $n\mid X\mid$.

After the  
spanning tree algorithm is run, 
the list of the points of each 
tree is intersected with itself
to avoid eventual repetitions that are left. 
At the end, each connected
component is 
associated to its tree in exactly one way.

In rare cases, the splitting step may return only one tree, but these cases 
are detected and repaired easily.

\subsection{Phase 2: Blow-up}\RSlabel{SecP1BlU}
This 
is also an adaptive iterative process. It reduces the set
$
\displaystyle{ X^u:=X\setminus \cup_{j=1}^J X_j^g}
$ 
of ``unsure'' data points
gradually, moving points from $X^u$ to one of the nearest sets $X_j^g$.
In order to deal with easy cases first, the points $\bx_i$ in $X^u$ 
are sorted by their locality quality 
such that points with
better localization come first. We also assume that for each point
$\bx_i\in X^u$ we know its distance to all sets $X_j^g$, and we 
shall update this distance during 
the algorithm, when the sets $X^u$ and $X_j^g$
change. We also use 
the distances to the sets
 $X_j^{g,0}$ that are the output of the localization phase and serve as a
  start-up for the sets $X_j^g$.

In an outer loop we run over all points $\bx_i\in X^u$ with 
decreasing quality of local approximation. In our implementation, this means 
increasing values of $\sigma_i$.
The inner loop runs over the $m$ sets $X_j^g$ to which
$\bx_i$ has shortest distance. In most cases, 
and in particular in $\R^2$, it will suffice to take $m=2$. 
The basic idea is to find the nearby set $X_j^g$ 
of ``good'' points for which the addition of $\bx_i$ 
does least damage to the local approximation quality.

Our implementation of the inner loop over $m$ 
neighboring sets $X_j^g$ 
works as follows.
In $X_j^{g,0}$, the point $\boy_j$ with shortest 
distance to $\bx_i$ is picked,
and its $n$ nearest neighbors in $X_j^{g,0}$ 
are taken, forming a set $Y_j^g$.
On this set, the data interpolant $s_j^g$ is calculated, and 
then the number $\sigma_j^g:=\|s_j^g\|_K$ 
measures the local approximation quality
near the point $\boy_j$ if only ``good'' points are used.
Then the ``unsure'' point $\bx_i$ is taken into account by forming a set 
$Y_j^u$ of points consisting of $\bx_i$  and the up to $n-1$ nearest neigbors to
$\bx_i$ from $X_j^g$. On this set, the data interpolant $s_j^u$ 
is calculated, and the  number $\sigma_j^u:=\|s_j^u\|_K$ 
measures the local approximation quality
if the ``unsure'' point $\bx_i$ is added to $X_j^g$. 
The inner loop ends by maintaining the minimum of quotients
$\sigma_j^u/\sigma_j^g$ over all nearby sets $X_j^g$ checked by the loop.
These quotients are used to indicate how much the local approximation quality
would degrade if $\bx_i$ would be added to $X_j^g$. Note that 
this strategy maintains locality by focusing on ``good'' nearest 
neighbors of either $\bx_i$ or $\boy_j$. By using the fixed sets 
$X_j^{g,0}$ instead of the growing sets
$X_j^g$, the algorithm does not rely heavily on the newly added points. 

An illustration is attached to Example $1$ in the next section; there
the point
$\boy_1$ and the sets  $X_1^g$ and $Y_1^u$
associated to a point $\bx_i$ will be shown.

After the inner loop, if 
the closest set to  ${\bf x}_i$ 
among all sets $X_k^g$ is $X_j^g$ and $\sigma_j^u/\sigma_j^g$
is less than $\sigma_k^u/\sigma_k^g$ for $k\neq j$,
then ${\bf x}_i$ is moved  from $X^u$ to $X_j^g$. If  
the closest set to  ${\bf x}_i$ is $X_j^g$ but if it is not true that
$\sigma_j^u/\sigma_j^g$ is less than $\sigma_k^u/\sigma_k^g$ for $k\neq j$,
then ${\bf x}_i$ remains ``unsure''.
The  ``unsure'' points are those that seriously degrade the
local approximation quality of all nearby sets of ``good'' points.

\subsection{Phase $3$: Final Assignment}

\RSlabel{SecP1FV}

The assignment of a point ${\bf x}_i\in X^u$ 
to a set $X_j^g$ is done on the basis
of how well the function value $f({\bf x}_i)$ is 
predicted by $u_j({\bf x}_i)$. We loop over all points $\bx_i\in X^u$
and first determine two sets $X_j^g$ and $X_k^g$
to which ${\bf x}_i$ has shortest distance. This is done in order to make sure
that $\bx_i$ is not assigned to a far-away $X_j^g$.
We then could add ${\bf x}_i$  
to $X_j^g$ if 
$|f({\bf x}_i)-u_j({\bf x}_i)|\leq |f({\bf x}_i)-u_k({\bf x}_i)|$,
otherwise to $X_k^g$, but in case that we have more than one unsure point,
we want to make sure that 
under all unsure points, $\bx_i$ fits better into $X_j^g$ than into $X_k^g$.
Therefore we calculate
$$
\begin{array}{rcl}
d_j(\bx_i)&:=& |f({\bf x}_i)-u_j({\bf x}_i)|\\
\mu_j&:=&\displaystyle{\min_{\bx_i\in X^u} d_j(\bx_i)   }\\ 
M_j&:=&\displaystyle{\max_{\bx_i\in X^u} d_j(\bx_i)   }\\ 
D_j(\bx_i)&:=&\displaystyle{\frac{d_j(\bx_i)-\mu_j}{M_j-\mu_j}   } 
\end{array}
$$
for all $j$ and $i$ beforehand, and assign ${\bf x}_i$  
to $X_j^g$ if $D_j(\bx_i)\leq D_k(\bx_i)$, otherwise to $X_k^g$.

\section{Examples}\RSlabel{SecExamples}

Some test functions are considered now, each of which is smooth 
on $J=2$ subdomains of $\Omega$. The algorithm constructs 
$X_1^g$ and $X_2^g$ with 
$X_1^g\cup X_2^g=X$.

Concerning  
the error of approximation of $u$, we separate what happens away
from the boundaries of  $\Omega_j$ from what happens globally on $[0,1]^2$.
This is due to the fact that standard domain boundaries, even without any
domain splittings, let the  approximation quality decrease near the boundaries.

To be more precise,  
let $\Omega_{safe}$ be the union of the circles of radius 
$$
q:=\displaystyle{ \min_{1\leq i<j\leq N}\|\bx_i-\bx_j\|_2},
$$
the separation distance of the data sites,
centered at those points of $X_j^g,\; j=1,2$ such that the 
centered circles of radius $2 {q}$ do not contain points 
of $X_k^g$ with $k\neq j$. 
We then  
evaluate
\[
L_\infty^{safe}(u):=\| u-f\|_{\infty,\Omega_{safe}\cap 
[0,1]^2}
\] and
\[
L_\infty(u):=\| u-f\|_{\infty,[0,1]^2}.
\]
 The chosen kernel for calculating 
the local kernel-based interpolants
is  the inverse multiquadric
 kernel $\phi(r)=(1+2 r^2/{\delta^2})^{-1/2}$ with parameter $\delta=0.35$.

 In all cases,  $N=900$ data locations are mildly scattered 
 on a domain $\Omega$ that extends 
 $[0,1]^2$ a little, with $q=0.04$. We shall restrict to $[0,1]^2$
 to evaluate the subapproximant, calculated by 
 the basis in the 
 Newton form. Such a basis is much more stable than the 
standard basis, see \RScite{mueller-schaback:2009-1}. 
The  error is computed on a grid with   step 
$0.01$. 
\par\bigskip\noindent{\bf Example 1}.
The function   
\[f_1(x,y):= \log(\mid x-(0.2 \sin(2\pi y)+0.5)\mid+0.5),\]
 has a  
derivative discontinuity 
across  the curve  $x=0.2 \sin(2\pi y)+0.5$. We get
$$
L_\infty^{safe}(u)=1.6\cdot  10^{-5},\;\;
L_\infty(u)=6.0\cdot 10^{-2}.
$$

For comparison, the errors of the global interpolant
are
$$
L_\infty^{safe}(u^{\star})=1.1\cdot  10^{-1},\;\; 
L_\infty(u^{\star})=1.1\cdot 10^{-1}.
$$

The classification turns out to be  
correct.   
$890$ out of  $900$ data  points are 
correctly classified as output of phase $3.2$, and then phase $3.3$ 
completes the classification.

Figure \ref{f1new} shows the points 
of $X_1^f$ as dotted and those of $X_2^f$ as crossed. 
The points both dotted and circled
of $X_1^f$, respectively the points both crossed and circled of $X_2^f$, 
are the result of the
splitting (Section \ref{SecP1S}), while the points dotted only,
respectively crossed only, are those
added by the blow-up phase (Section \ref{SecP1BlU}). The points squared
are the result of the final assignment
phase (Section \ref{SecP1FV}).
The true splitting line is traced too. The convention of 
the marker types will be used 
in the next examples as well.

The function $u$ is defined as $u_1^f$ where 
the subdomain $\Omega_1$ is determined
and as $u_2^f$ on $\Omega_2$. 

The actual error $L_\infty(u)=6.0\cdot 10^{-2}$ is not much affected if we omit
Phase 3 and and ignore the remaining 10 ``unsure'' points after the blow-up
phase. 
 A similar effect is observed for the other
  examples to follow.

A zoomed area of $\Omega$ is considered in 
Figure  \ref{f4zoom}. The details are related to an
iteration of the blow-up phase, 
where  the ``unsure'' point ${\bf x}_i$ (both 
squared and starred)  is currently
examined.
Points of
${X}_2^{g,0}$ are shown  as crosses.
At the current iteration,  the  dots are 
points inserted in 
${X}_1^{g}$ up to now, those belonging to ${X}_1^{g,0}$ bold dotted, 
while the points inserted in ${X}_2^g$ up to now
are omitted in this illustration.
The 
points squared  are those of ${ Y}_1^u$, while the points as
diamonds
are those of ${ Y}_1^g$. The point ${\bf y}_1$ is both 
written 
as diamond and
star. 


\par\bigskip\noindent{\bf Example 2}.
The function
\[\RSlabel{fun1modif}
f_2(x,y):=\left\{\begin{array}{ll}
f_1(x,y)  \quad\quad\quad {\rm if}\quad x<=0.2\, \sin(2\pi y)+0.5 \\
f_1(x,y)+0.01 \quad {\rm if}\quad x>0.2\, \sin(2\pi y)+0.5
\end{array}
\right.\]

has a discontinuity across the curve $x=0.2\,\sin(2\pi y) + 0.5$.

We get
\[L_\infty^{safe}(u)=1.6\cdot 10^{-5},\;L_\infty(u)=6.0\cdot 10^{-2}.\]

For comparison, the errors of the global interpolant are
\[L_\infty^{safe}(u^\star)=1.3\cdot 10^{-1}
,\;L_\infty(u^\star)=1.3\cdot 10^{-1}.\]
The classification turns out to be correct.
 $888$ out of $900$ data points are classified correctly
as output of phase $3.2$, and phase $3.3$ completes the classification
for the remaining 12 points.
It might be that $u_j^f$ is more accurate on the safe zone, and also globally.

\par\bigskip\noindent{\bf Example 3}.
The function  
\begin{equation}\RSlabel{fun3}
f_3(x,y):= \arctan(10^3 (\sqrt{(x+0.05)^2+(y+0.05)^2}-0.7)) 
\end{equation}
is regular but  has a steep gradient.  
Our algorithm yields
$$
L_\infty^{safe}(u)=9.0\cdot 10^{-2} \hbox{ and } L_\infty(u)= 2.67\cdot 10^0,
$$
while for the global interpolant we get 
\[L_\infty^{safe}(u^\star)=2.31\cdot 10^{0}
,\;L_\infty(u^\star)=3.26\cdot 10^{0}.\]
 Figure \ref{sigm1} shows the points of $X_1^f$
as dotted and those of $X_2^f$ as crossed; $X_1^f$ and $X_2^f$
stay at the opposite sides of the mid range line
$f(x,y)=0$ .

\par\bigskip\noindent{\bf Example 4}.
The function 
\begin{equation}\RSlabel{fun4}
f_4(x,y):= ((x-0.5)^2+(y-0.5)^2)^{0.35}+0.05*(x-0.5)^0_+ 
\end{equation}
has a jump on the line $x=0.5$ and a derivative singularity on it at
$(0.5, 0.5)$.  It has rather a steep gradient too.
One data point close to the singularity is not
classified correctly. We get
$$
L_\infty^{safe}(u)=9.9\cdot 10^{-4} \hbox{ and } L_\infty(u)=7.3\cdot 10^{-2},
$$
while the global interpolant $u^{\star}$ has 
$$
L_\infty^{safe}(u^{\star})=1.5\cdot 10^{-2} \hbox{ and }  
L_\infty(u^{\star})= 8.5\cdot 10^{-2}.
$$
 Figure \ref{punta1} shows the points of $X^f_1$ as dotted
and those of $X^g_2$ as crossed.

All examples show that the transition from a global to a
properly segmented problem decreases the achievable error considerably.
But the computational cost is serious, and it might be more efficient 
to implement a multiscale strategy that works on coarse data first,
does the splitting of the domain  coarsely, and then refines
the solution on more data, without recalculating everything on the finer data. 

\begin{figure}[hbt]
\psfig{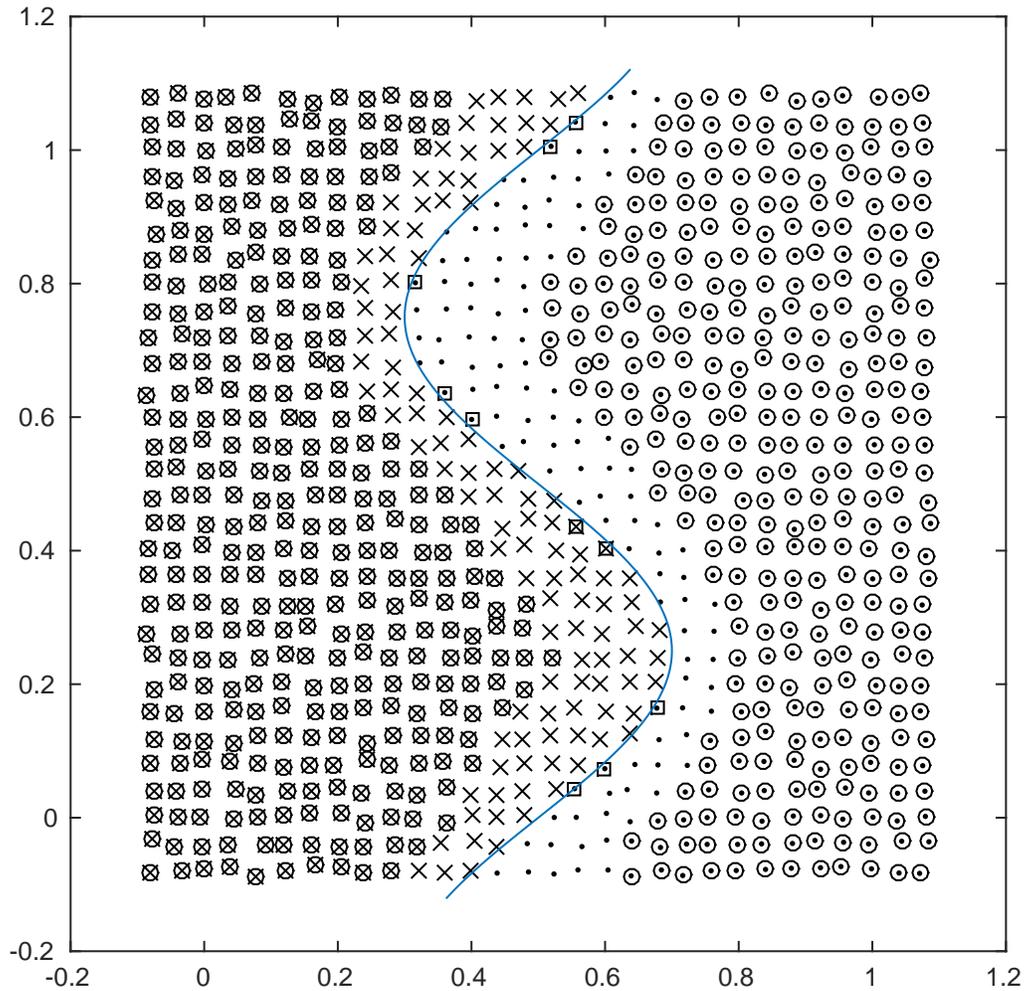}\\
\caption{Example $1$: class $1$ as dots, class $2$ as crosses}
\RSlabel{f1new}
\end{figure}
\begin{figure}[hbt]
\psfig{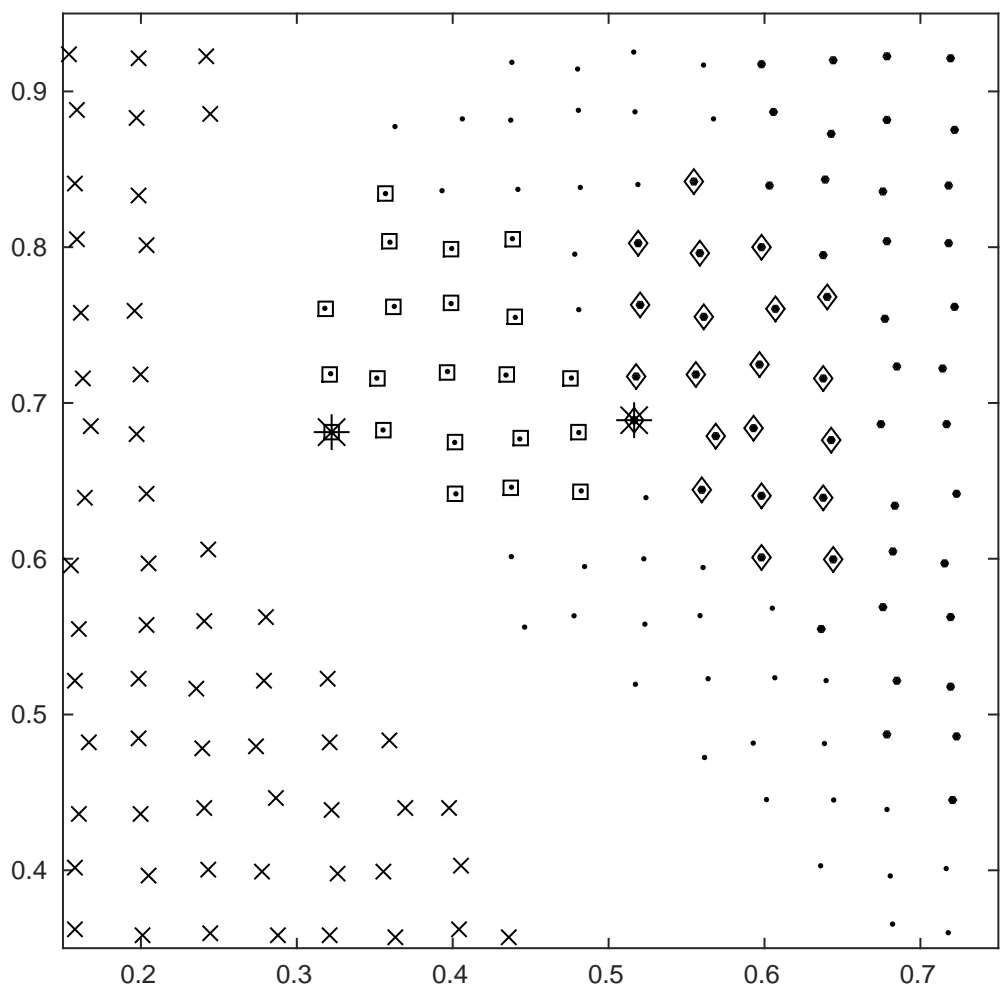}\\
\caption{Localized blow-up, zoomed in }
\RSlabel{f4zoom}
\end{figure}
\begin{figure}[hbt]
\psfig{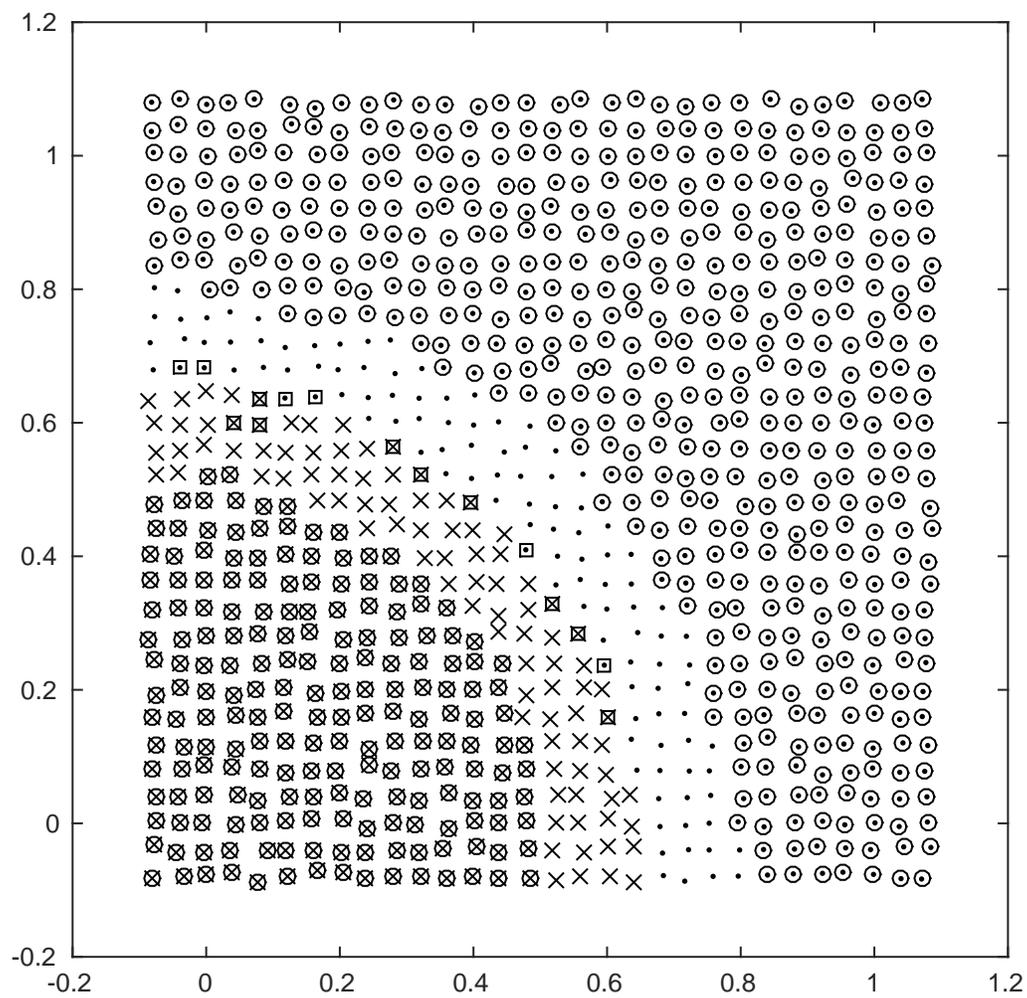}\\
\caption{Example $3$:  class $1$ as dots, class $2$ as crosses}
\RSlabel{sigm1}
\end{figure}
\begin{figure}[hbt]
\psfig{figure=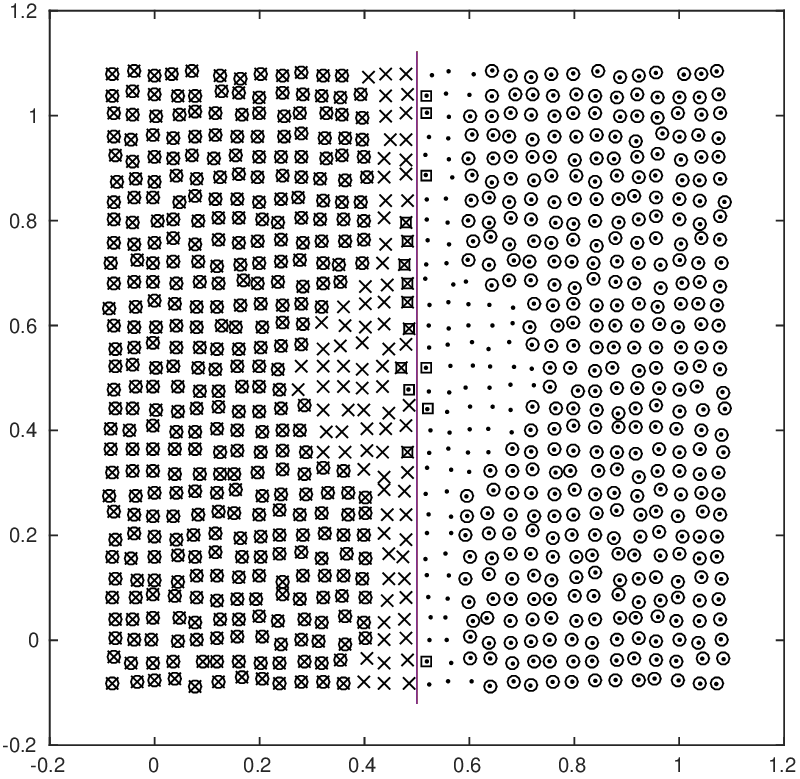}\\
\caption{Example $4$:   class $1$ as dots, class $2$ as crosses}
\RSlabel{punta1}
\end{figure}

\bibliographystyle{plain}

\end{document}